\documentclass[10pt]{article} 
\usepackage{amsfonts,amssymb,latexsym,amsmath,oldgerm,amscd}

\newcommand{\mf}{\mathfrak}

\newcommand{\ra}{\rightarrow}
\newcommand{\Ra}{\Rightarrow}

\newcommand{\eps}{\varepsilon}

\newcommand{\mc}{\mathcal}
\newcommand{\mbb}{\mathbb}
\newcommand{\tn}{\textnormal}

\newtheorem{de}{Definition}[section]

\newtheorem{tr}[de]{Theorem}
\newtheorem{lm}[de]{Lemma}

\newtheorem{co}[de]{Corollary}

\def\vp{\rm \vspace{0.1cm}}
\def\M{\rm Max}
\def\hb{\hfill$\Box$}
\def\ZZ{\mathbb{Z}} 
\def\M{\rm Max}
\def\m{\mf{m}}

\def\hb{\hfill$\Box$}

\def\GL{\rm GL}
\def\SL{\rm SL}

\def\E{\rm E}

\def\t{\rm K_1}
\def\K{\rm K}

\def\nk{\rm NK_1}
\def\k0{\rm NK_0}
\def\w{\rm W}
\def\sk{\rm SK_1}

\def\m{\rm M}

\title{Absence of torsion for ${\nk}(R)$ over associative rings}
\author{Rabeya Basu \\{\small {\it Department of Mathematical Sciences,}}\\
{\small {\it Indian Institute of Science Education and Research, Kolkata,}} \\{\small {\it West Bengal, India}}\\
{\small {\it rbasu@iiserkol.ac.in}}}

\begin{document}

\date{}
\maketitle
\begin{center}{\it 2000 Mathematics Subject Classification:{13H99, 15A24, 16R50, 19B14}}
\end{center}
\begin{center}
{\it Key words: linear group, ${\t}$, ${\nk}$, ${\sk}$, torsion, Witt vectors.}
\end{center}
\begin{center} {\bf Abstract} 
\end{center}
When $R$ is a commutative ring with identity, and if $k \in \mbb{N}$, with $kR = R$, then it 
was shown in \cite{WEL} that ${\sk}(R[X])$ has no k-torsion. We reprove this
result for any {\it associative} ring $R$ with identity in which $kR = R$.


\section{Introduction} 
~~~~Let $R$ be an associative ring with identity element $1$. 
Let ${\t}(R)$ denote the Whitehead group. In case $R$ is commutative, let ${\sk}(R)$ be the kernel of 
the determinant map from ${\t}(R)$ to the group of units of $R$. 
Let ${\w}(R)$ be the ring of big Witt vectors. We denote ${\nk}(R)=\tn{ker}
({\t}(R[X])\ra {\t}(R))$; $X=0$.
In \cite{STI}, J. Stienstra, using ideas of S. Bloch in \cite{BL1}, showed that 
${\nk}(R)$ is a ${\w}(R)$-module. Consequently, as noted by C. Weibel in 
(\cite{WEL}, \S 3), if $k$ is a unit in $R$, then ${\sk}(R[X])$ has no 
$k$-torsion, when $R$ is a commutative local ring. For a commutative ring with identity
${\nk}(R)$ coincides with
${\sk}(R[X])$ if we take $R$ to be local. 
In this note we generalize Weibel's observation for ${\nk}(R)$, where $R$ is a commutative local ring.

We prove for an \textit{associative} ring $R$ with identity, ${\nk}(R)$ has no $k$-torsion if $k$ is a unit in $R$. In particular, this shows Weibel's result is a special case. 
The method of proof may be considered as a simplified version of J. Stienstra's approach via big Witt vectors. This will help the reader to appreciate why big 
Witt vectors come into the picture in the functorial approach of S. Bloch, 
{\it et al}. \vp \\
\begin{tr}  \label{th1} Let $R$ be an associative ring with identity. If $k$ is a unit 
in $R$, then ${\nk}(R)$ has no $k$-torsion.
\end{tr}
\begin{co}
Let $R$ be a commutative local ring with identity. If $k$ is a unit in $R$, then
${\sk}(R[X])$ has no $k$-torsion.
\end{co}

As a consequence of Theorem 1 we prove 
\begin{tr} \label{th2} 
Let $R=R_0\oplus R_1\oplus \cdots $ be a graded commutative ring with identity. 
Let $k$ be a unit in $R_0$.  
Let $N=N_0+N_1+\cdots+N_r \in {\m}_r(R)$ be a nilpotent matrix.
If $[(I+N)]^k=[I]$ in ${\sk}(R)$, then $[I+N] = [I+N_0]$. 
In particular, if $R_0$ is a reduced local ring, then ${\sk}(R)$ has no 
$k$-torsion. 
\end{tr}

\section{Prologue} 
~~ Let $R$ be an associative ring with identity.
${\GL}_n(R)$ denotes the group of invertible matrices, ${\SL}_n(R)$ its subgroup 
of matrices of determinant $1$ (when $R$ is a commutative ring), 
${\E}_n(R)$ the subgroup of elementary matrices, 
{\it i.e.} generated by $\{{\E}_{ij}(\lambda):\lambda \in R, i\ne j \}$, where 
${\E}_{ij}(\lambda)=I+\lambda e_{ij}$ and $e_{ij}$ is the matrix with $1$ on the 
$ij$-th position and 0's elsewhere. For $\alpha\in {\m}_r(R)$, $\beta\in {\m}_s(R)$ 
we have $\alpha\perp \beta\in {\m}_{r+s}(R)$, where $$\alpha\perp \beta =
\left(\begin{array}{cc} \alpha & 0 \\ 0 & \beta
\end{array}\right).$$

There is an infinite counterpart: identifying each matrix $\alpha\in {\GL}_n(R)$ 
with the large matrix $(\alpha\perp 1)$
gives an embedding of ${\GL}_n(R)$ into ${\GL}_{n+1}(R)$. 
Let ${\GL}(R)=\cup_{n=1}^{\infty} {\GL}_n(R)$, 
${\SL}(R)=\cup_{n=1}^{\infty} {\SL}_n(R)$, and 
${\E}(R)=\cup_{n=1}^{\infty} {\E}_n(R)$ be the corresponding 
infinite linear groups.

The well known Whitehead's Lemma asserts that if $\alpha\in {\GL}_n(R)$ then we have 
$(\alpha\perp \alpha^{-1})\in {\E}_{2n}(R)$. Thus we have $$[{\GL}(R),{\GL}(R)]=[{\E}(R),{\E}(R)]={\E}(R)$$ and hence
${\E}(R)$ is a normal subgroup of ${\GL}(R)$. 
The quotient ${\GL}(R)/{\E}(R)$ is called the
{\bf Whitehead group} of the ring $R$ and is denoted by ${\t}(R)$. 
For $\alpha\in {\GL}_n(R)$ let $[\alpha]$ denote its equivalence class in 
${\t}(R)$. 
Also, as a consequence of Whitehead's lemma
one sees that if $\alpha,\beta\in {\GL}_n(R)$ then 
$[\alpha,\beta]\in {\E}_{2n}(R)$;
whence ${\t}(R)$ is an abelian group. For details
{\it cf.} \cite{B}.

In case $R$ is commutative the determinant map from 
${\GL}_n(R)$ to $R^{*}$ induces a map,
$\det: {\t}(R) \ra R^{*}$ given by $\alpha E(R)\mapsto \det \alpha$. 
The kernel of the map is denoted by ${\sk}(R)$ and equals to ${\SL}(R)/{\E}(R)$.

We  write ${\nk}(R)$ for $\tn{ker} ({\t}(R[X])\ra {\t}(R))$; $X=0$,  
{\it i.e.} the subgroup consisting of elements $[\alpha(X)]\in {\t}(R[X])$ 
such that $[\alpha(0)]=[I]$. Note that if $R$ is a commutative local ring then 
${\sk}(R[X])$ coincides with ${\nk}(R)$; indeed, if $R$ is a local ring then
${\SL}_n(R)={\E}_n(R)$ for all $n>0$. Therefore, we may replace $\alpha(X)$ by 
$\alpha(X)\alpha(0)^{-1}$ and assume that $[\alpha(0)]=[I]$. \vp\\

For a commutative ring $R$ the group ${\w}(R)$ of big {\bf Witt vectors} is 
defined by: $${\w}(R)=(1+XR[[X]])^{\times}.$$ 
For $P(X)\in (1+XR[[X]])$, let $\omega(P)$ denote the corresponding 
element of ${\w}(R)$. The group structure in ${\w}(R)$ is given by: 
$$\omega(P)+\omega(Q)=\omega(P.Q).$$ 
Any $P(X)\in (1+XR[[X]])$ can be written uniquely as a 
product: $$P(X)=\underset{n\ge 1}\Pi (1-a_nX^n)^{-1}, \,\, a_n\in R.$$
The elements $(a_1,a_2,\dots \dots)$ are called the 
{\bf Witt co-ordinates} of $\omega(P)$. Also, there exists a unique structure 
of commutative ring on ${\w}(R)$ such that 
$$\omega ((1-aX^m)^{-1}).\omega ((1-bX^n)^{-1})=
\omega ((1-a^{n/r}b^{m/r}X^{mn/r})^{-r}),$$
where $r=\tn{ g.c.d} (m,n)$. The identity element in ${\w}(R)$ is presented by the power 
series $(1-X)^{-1}$. For details {\it cf.} (\cite{BL1}, Prop. I.I), \cite{L}. 

For the ${\w}(R)$-module structure of ${\nk}(R)$ see \cite{WEL1}
(for previous articles {\it cf.} \cite{BL1}, \cite{BL2}, \cite{STI},
\cite{WEL}).
\section{Higman Linearization} 
~~Two matrices $\alpha\in {\m}_r(R)$ and $\beta\in {\m}_s(R)$ are said 
to be {\bf stably equivalent} if there exists $\eps_1, \eps_2\in {\E}_t(R)$ 
(for some $t\ge \tn{max} \{r,s\}$) such that 
$\eps_1 (\alpha \perp I_{t-r})\eps_2 = (\beta \perp I_{t-s})$.
\begin{lm} \label{hig} {\bf (Higman Linearization Process)} 
Let $\alpha(X)$ be a matrix over $R[X]$. Then $\alpha(X)$ is 
stably equivalent to a linear matrix in ${\m}_s(R[X])$ for some $s$. 
\end{lm} 
{\bf Proof.} We may assume that $n\ge 2$. 
Let $$\alpha(X)=a_0+a_1X+a_2X^2+\cdots +a_nX^n\in {\m}_r(R), \,\,\, r>1.$$ 
Then $\alpha(X)$ is stably equivalent to a matrix of degree $n-1$ over 
$R[X]$ in the following manner:
$$\left(\begin{array}{cc}
I_r & -a_nX\\
0  & I_r
\end{array} \right)
\left(\begin{array}{cc}
a_0+\cdots +a_nX^n & 0\\
0 & I_r 
\end{array} \right)
\left(\begin{array}{cc}  
I_r & 0\\
X^{n-1}I_r & I_r
\end{array} \right) \\ $$ 
$$=
\left(\begin{array}{cc}
a_0+\cdots +a_{n-1}X^{n-1} & -a_nX\\
X^{n-1}I_r & I_r 
\end{array} \right)=\alpha_1$$ 
has degree $(n-1)$. Hence $\alpha$ is stably equivalent to $\alpha_1$.
Repeating the above process $(n-2)$ times we get the result. \hb

\begin{co}  \label{hig1} Let $R$ be an associative ring with identity. 
Let $\alpha(X)\in {\GL}_r(R[X])$ with $\alpha(0)=I_n$. Then in ${\t}(R[X])$ 
we have $[\alpha(X)]=[I_s+NX]$ for some $s>0$ and some matrix $N\in \m_s(R)$.  
\end{co} 
{\bf Proof.} By Lemma \ref{hig} 
there exists $\eps_1, \eps_2\in {\E}_t(R[X])$ (for some $t>r$) such that 
$(\alpha(X) \perp I_{t-r})= \eps_1((I_s+NX) \perp I_{t-s})\eps_2$ for 
some $s>0$ and $N\in \m_s(R)$. Now as ${\E}(R)$ is a normal subgroup of 
${\GL}(R)$,  for some integer $u$ there exists $\eps_1'\in {\E}_{t+u}(R[X])$ such that 
$$(\eps_1\perp I_u)((I_s+NX) \perp I_{t-s+u})= 
((I_s+NX) \perp I_{t-s+u})\eps_1'.$$ 
Hence in ${\t}(R[X])$ we have 
$$[\alpha(X)]=[((I_s+NX) \perp I_{t-s+u})\eps_1'(\eps_2\perp I_u)]=
[I_s+NX].$$ \hb

\section{Main Theorem} 

~~Let $R_t$ denote the ring $R[X]/(X^{t+1})$.

\begin{lm} \label{la3}  
Let $R$ be a ring and $P(X)\in R[X]$ be any polynomial.
Then the following identity holds in the ring $R_t:$ 
\begin{equation*}
(1+X^r P(X))=(1+X^rP(0))(1+X^{r+1}Q(X)),
\end{equation*}
where $r>0$ and  $Q(X)\in R[X]$, with $\deg(Q(X))< t-r$.
\end{lm}
 {\bf Proof.} Let us write $P(X)=a_0+a_1X+\cdots+a_{t}X^{t}$. Then we can 
write $P(X)=P(0)+XP'(X)$ for some $P'(X)\in R[X]$. Now, in $R_t$
{\small
\begin{align*}
(1+X^r P(X))(1+X^r P(0))^{-1} 
& =  (1+X^r P(0)+X^{r+1}P'(X))(1+X^r P(0))^{-1}\\
& = 1+X^{r+1}P'(X)(1-X^rP(0)+X^{2r}(P(0))^2-\cdots)\\
& =  1+X^{r+1}Q(X)
\end{align*}}
\!\!where $Q(X)\in R[X]$ with $\deg(Q(X))< t-r$. 
Hence the lemma follows.  \hfill$\Box$ \vp \\
{\bf Remark.}  Iterating the above process we can write for any polynomial 
$P(X)\in R[X]$,
$(1+XP(X))=\Pi_{i=1}^t(1+a_iX^i)$ in $R_t$, 
for some $a_i\in R$. By ascending induction it will follow that the $a_i$'s 
are uniquely determined. In fact, if $R$ is commutative then 
$a_i$'s are the $i$-th component of the 
ghost vector corresponding to the big Witt vector of 
$(1+XP(X))\in {\w}(R)=(1+XR[[X]])^{\times}$.  For details see 
(\cite{BL1}, $\mc{x}$I).

 \begin{lm} \label{la4} 
Let $R$ be a ring with $\frac{1}{k}\in R$ and  $P(X)\in R[X]$. 
Assume $P(0)$ lies in the center of $R$. 
Then 
$$(1+X^rP(X))^{k^r}=1 \Ra (1+X^rP(X))=(1+X^{r+1}Q(X))$$ in the ring $R_t$
for some $r>0$ and $Q(X)\in R[X]$ with $\deg(Q(X))< t-r$.
\end{lm} 
{\bf Proof.} By Lemma \ref{la3} 
\begin{equation} \label{use4} 
(1+X^rP(X))=(1+X^rP(0))(1+X^{r+1}P_1(X)),
\end{equation}
for some $P_1(X)\in R[X]$ with $\deg(P_1(X))< t-r$. Therefore, in $R_t$
$$(1+X^rP(X))^{k^r}=1\Ra   
(1+X^rP(0))^{k^r}=(1+X^{r+1}P_1(X))^{-{k^r}}.$$
As $\frac{1}{k}\in R$, we have 
$$(1+{k^r}X^rP(0)+X^{r+1}P_2(X))=(1+X^{r+1}P_1(X))^{-{k^r}}.$$ This implies 
\begin{align*}
(1+{k^r}X^rP(0))
  & =  \,\, (1+X^{r+1}P_1(X))^{-{k^r}} 
(1+(1+{k^r}X^rP(0))^{-1}X^{r+1}P_2(X))^{-1}\\
& = \,\, (1+X^{r+1}P_3(X))
\end{align*}
for some $P_2(X),P_3(X)\in R[X]$ with deg $P_2(X),P_3(X)< t-r$. Now,
applying homomorphism $X\mapsto \frac{1}{k}X$ we get
$$(1+X^rP(0))=(1+X^{r+1}P_4(X))$$ for some $P_4(X)\in R[X]$ 
with $\deg(P_4(X))< t-r$. 
Substituting this in \eqref{use4} we get 
$$(1+X^rP(X))=(1+X^{r+1}Q(X))$$ for some $Q(X)\in R[X]$ with 
$\deg (Q(X))< t-r$. \hb \vp \\
{\bf Proof of Theorem} \ref{th1}.
Let $\alpha(X)\in {\GL}_n(R[X])$ with $[\alpha(0)]=[I]$ be a $k$-torsion. 
By Corollary \ref{hig1} in ${\t}(R[X])$,
$[\alpha(X)]=[(I_s+NX)]$ for some $s>0$ and 
$N\in {\m}_s(R[X])$. 
Since $(I_s+NX)$ is invertible, $N$ is nilpotent.  
Let $N^{t+1}=0$. Since $[(I_s+NX)]^k=[I]$ in ${\t}(R[X])$, it follows that 
$$[I_s+kNX+N^2X^2P_1(NX)] = [I]$$ for some $P_1(X)\in R[X]$. Hence, as before, 
as $\frac{1}{k}\in R$,
\begin{align*}
[(I_s+kNX)^{-1}] & = [(I_s+(I_s+kNX)^{-1}N^2X^2P_1(NX)]\\ 
& = [I_s+(I_s-kNX+N^2X^2P_2(NX))N^2X^2P_1(NX)]\\
& = [I_s+N^2X^2P(NX)]
\end{align*}
for some $P(X)\in R[X]$. 
Since $[(I_s+N^2X^2P(NX))]^k=[I]$, arguing as in the proof of Lemma \ref{la4} 
we get in ${\t}(R[X])$ 
$$[I_s+N^2X^2P(NX)]=[I_s+N^3X^3Q(NX)]$$ for some $Q(X)\in R[X]$. 
Now by repeating the above mentioned 
argument  we get $$[I_s+N^2X^2P(NX)]=[I].$$ 
Finally, applying homomorphism 
$X\mapsto \frac{1}{k}X$ we get the desired result.\hb \vp \\

Now we prove Theorem \ref{th2} as a consequence of Swan-Weibel homotopy trick. 
For details see (\cite{GUB1}, Proof of Prop. 2.22).
First we recall the Local-Global Principle for a graded ring. 
\vp \\ 
{\bf Graded Local-Global Principle:}
Let $R=R_0\oplus R_1\oplus \cdots$ be a graded commutative ring with 
$k$ a unit in $R_0$ and $\alpha(X)\in {\GL}_n(R[X])$ with $\alpha(0)=I_n$. If  
$\alpha_{\mf{m}}(X)\in {\E}_n(R_{\mf{m}}[X])$ for all $\mf{m} \in {\M} \, (R_0)$, 
then $\alpha(X)\in {\E}_n(R[X])$. \vp

This is derived from the usual Local-Global Principle by using the following 
homotopy trick due to Swan and Weibel.
\vp \\ 
{\bf Proof of Theorem \ref{th2}.}
Consider the ring homomorphism  
$\theta: R\ra R[X]$, given by 
$(a_0+a_1+\cdots)\mapsto a_0+a_1X+\cdots$.
Then 
\begin{align*}
[(I+N)]^k=[I] &\Ra \theta([(I+N)]^k)=[(\theta(I+N))]^k=[I]\\ 
& \Ra [(I+N_0+N_1X+\cdots+N_rX^r)]^k=[I]
\end{align*}
Let $\mf{m}$ be a maximal ideal in $R_0$. 
By Theorem 1
\begin{equation*}[(I+N_0+N_1X+\cdots+N_rX^r)]=[I]
\end{equation*} in ${\sk}(R_0)_{\mf{m}}.$
Hence by the Graded Local-Global Principle $[(I+N)]=[(I+N_0)]$ in ${\sk}(R)$. 

In particular, if $R_0$ is a reduced local ring then 
units of $R_0[X]$ = units of $R_0$. 
Since $N_0$ is nilpotent, 
$\det(I+N_0X)$= constant = $1$. Hence 
$I+N_0$ is in ${\SL}_r(R_0)={\E}_r(R_0)$. This completes the proof.  \hb  

\section{Examples}
~~We show by an example that the the condition 
$\frac{1}{k}\in R$ is necessary in Theorem \ref{th2}, whence in 
Theorem \ref{th1}. \vp \\
{\bf Example.}
Let $R$ be a commutative ring with identity and  $\alpha$ a $2 \times 2$ 
completion of $(1 - XY, X^2)$ over $R[X^2, XY, Y^2]$. 
In  (\cite{GUB}, \S 8,  Example 8.2)  it is shown that $\alpha \in
{\sk}(R[X^2, XY, Y^2]) \setminus {\sk}(R)$. 
Now take $R = \ZZ_2$. 
Then the square of the Mennicke symbol of the vector 
$(1 - XY, X^2)$ is clearly trivial, 
{\it i.e.} $\alpha^2 \in {\E}_3(\mbb{Z}_2[X^2, XY, Y^2])$.
(See \cite{B} for definition of Mennicke symbol.) \vp \\
{\bf Remark.}
Let $P$ be a finitely generated projective $R$-module. A theorem of 
M.R. Gabel in \cite{GAB} asserts that  $mP=P\oplus \cdots \oplus P$ (m times) 
is free, for some $m$. (A result of T.Y. Lam in \cite{LAM} sharpens this bound 
on $m$). 

Ravi A. Rao has asked if the converse of the above is true over polynomial rings $R[X]$ 
with $R$ local. 
More precisely, if R is a local ring with kR = R, does ${\K}_0(R[X])$ have
non-trivial k-torsion?


\footnotesize
\setlength{\itemsep}{-0.5ex}
  
\addcontentsline{toc}{chapter}{Bibliography} \vp
{\it Indian Institute of Science Education and Research (IISER-K)},\\
{\it Mohanpur Campus, P.O. BCKV Campus Main Office},\\
{\it Mohanpur, Nadia - 741252, West Bengal, India.} \\
{\it Fax: 03473-2334-4107 }\\
{\it ~email : rabeya.basu@gmail.com, rbasu@iiserkol.ac.in}

\end{document}